# Green Vehicle Routing Problem: State of the Art and Future Directions

Saba Sabet[1a], Bilal Farooq[a]

[a]*Laboratory of Innovations in Transportation (LiTrans), Ryerson University*

**Abstract**

Green vehicle routing problem (GVRP) aims to consider greenhouse gas emissions reduction, while routing the vehicles. It can be either through adopting Alternative Fuel Vehicles (AFVs) or with existing conventional fossil fuel vehicles in fleets. GVRP also takes into account environmental sustainability in transportation and logistics. We critically review several variations and specializations of GVRP to address issues related to charging, pick-up, delivery, and energy consumption. Starting with the concepts and definitions of GVRP, we summarize the key elements and contributors to GVRP publications. Afterward, the issues regarding each category of green vehicle routing are reviewed, based on which key future research directions and challenges are suggested. It was observed that the main focus of previous publications is on the operational level routing decision and not the supply chain issues. The majority of publications used metaheuristic methods, while overlooking the emerging machine learning methods. We envision that in addition to machine learning, reinforcement learning, distributed systems, internet of vehicles (IoV), and new fuel technologies have a strong role in developing the GVRP research further.

## 1. Introduction

The Vehicle Routing Problem (VRP) is used to design optimal routes for a fleet of vehicles to service a set of customers, considering a certain set of constraints. VRP was first introduced by Dantzig and Ramser [1] in their seminal work on truck dispatching. In that study, the first algorithmic approach was proposed and applied to optimize fuel deliveries. According to the study, VRP aims to find optimal routes for a fleet of delivery trucks, each with limited capacity, so as to minimize the total distance travelled. There can be one or more depots and customer nodes in a vehicle routing network.

Green Vehicle Routing Problem (GVRP) is a branch of green logistics, which refers to vehicle routing problems where externalities of using vehicles, such as carbon dioxide-equivalents emissions, are to be reduced. In this way, although research in this field has a long history, Erdogan and Miller-Hooks [2] formally introduced the term GVRP for the first time. GVRP incorporates the environmental aspects of transportation into VRP, which is one of the most interesting problems in the field of logistics and transportation. The goal of this problem is to earn economic benefits, while also taking into account environmental considerations. In its most general form, GVRP aims to minimize GHG emissions with solely conventional gasoline and diesel vehicles (GDVs) or with Alternative Fuel Vehicles (AFVs) in the fleets. It also takes into account the environmental sustainability in freight transportation.

---

[1]Corresponding author

The aim of this study is to provide a systematic state of the art and outline new insights and perspectives into GVRP, based on a wide range of relevant searches by answering the main review questions below:

1. What are the main variants of GVRP developed to incorporate the environmental goals in the VRP field?
2. What are the strategies to address environmental issues in GVRP?
3. What multiple objectives associated with the GVRP variants should be considered?
4. What are the emerging issues and future trends in GVRP?

For the existing literature on GVRP, the availability of extensive resources from reputed journals, books, technical reports, surveys and conference proceedings helped the present study. Specifically, we collected the relevant reserach articles from academic databases, including Google Scholar, Scopus, Springer, Taylor & Francis Elsevier, ScienceDirect, Wiley, and IEEEXplore. To conduct the search, we identified keywords such as green vehicle routing, green logistics, electric vehicle routing,alternative fuel vehicle routing, and charging-discharging scheduling of electric vehicles.

The overall structure of the study is as follows: next section presents GVRP research background, Section 3 presents GVRP Algorithms, Section 4 presents GVRP classification, Section 5 presents new insights into GVRP, and the final section presents a detailed conclusion and future research direction.

## 2. Analysis of Existing Literature

The aim of this study is to provide a systematic literature review based on a wide range of relevant searches. The flowchart in Figure 1 refers to the review methodology applied in this article. The structure in the figure illustrates the current study in several steps. First, we determined the basis for the review work; retrieving and refining the selected papers and organizing the outline of the review paper. In the next step, we provided descriptive analysis using quantitative figures. Then, we used software and apply scientometric analysis on the refined papers to visualize the clustering of keywords used in the literature. In the next step, we focused on the fundamental aspects of GVRP, extracting various types of trends and the gap existing in the current GVRP research. In the end, the conclusion is described, as well as providing opportunities for future research on GVRP.

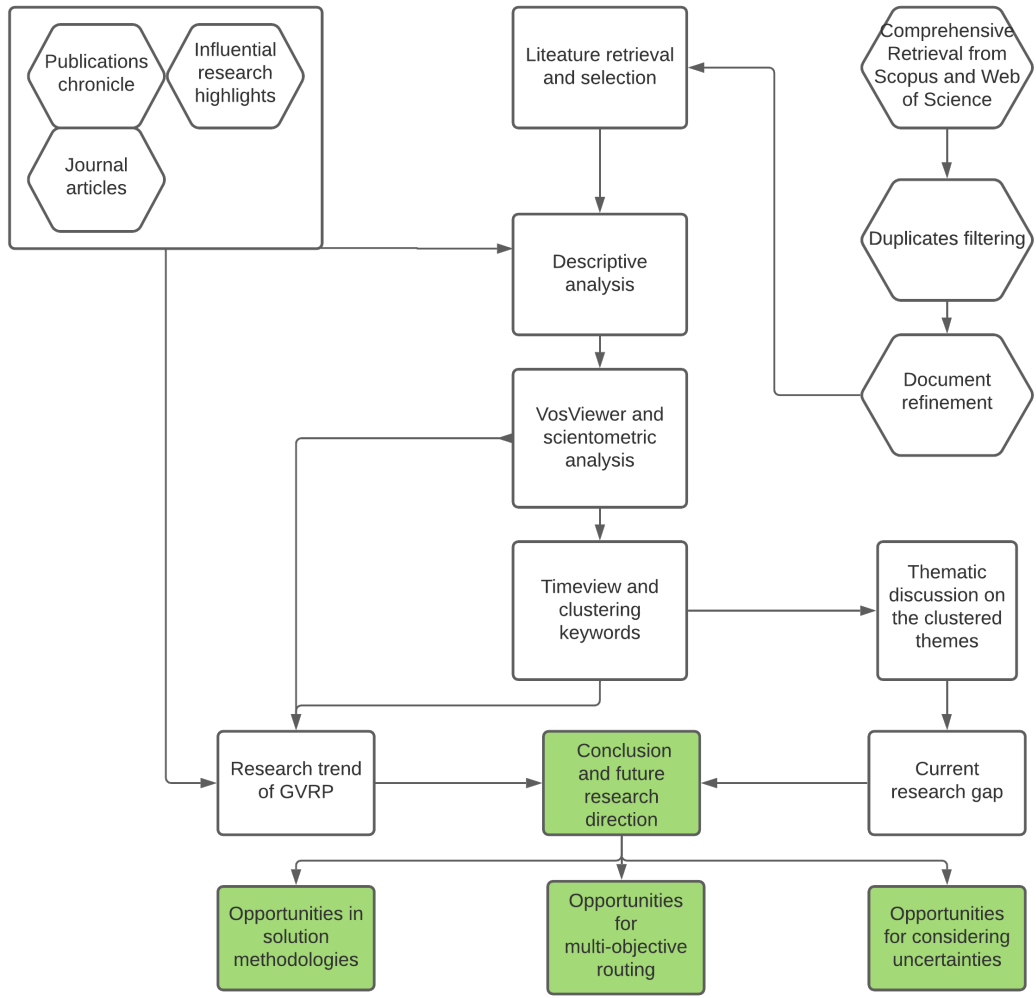

Figure 1: Procedure of applied literature analysis

## 2.1. Scientometric Analysis of GVPR Literature

Considering the importance of mitigating GHG emission, many studies have looked into green vehicle routing during the last decades using different terms such as Pollution Routing Problem (PRP), Eco-Routing and GVRP. Thus, it would be complicated to pinpoint a particular period in time when publications of green vehicle routing officially started. The earliest publications of this domain, as recorded by the web of science, date back to 2012 [2]. Since then, publications of GVRP have constantly been the center of attention with an overall number of more than 450 publications over the years. During the last five years leading to the time of the current publication, more than 400 publications on GVRP have been indexed each year, 100 items of which associated with 2019. In other words, the size of this literature has risen sharply within the last five years. See Figure 2 for a visualisation of this trend in GVRP literature.

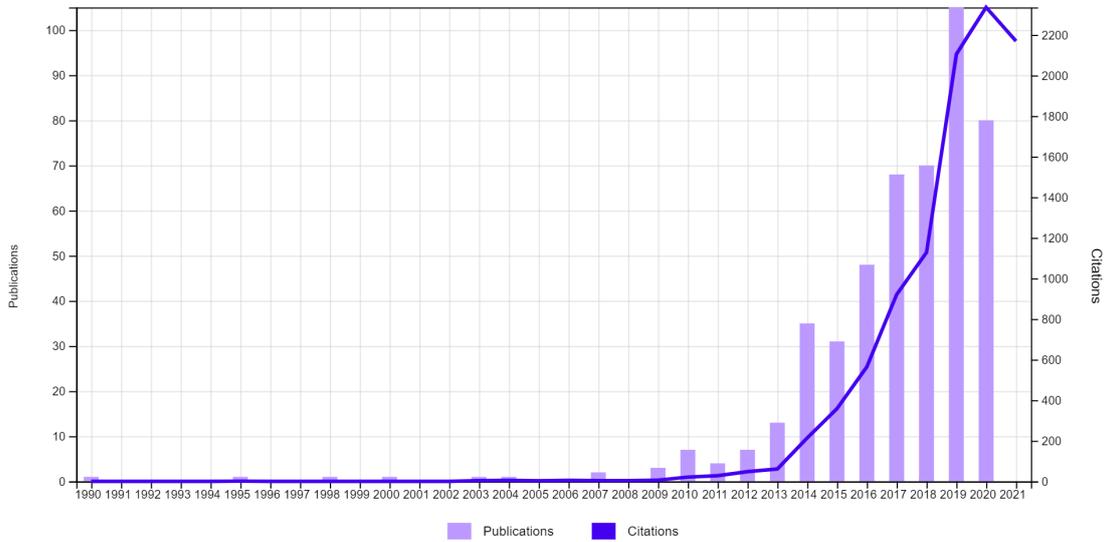

Figure 2: The number of GVRP publications by year

A significant number of studies are associated with scientific articles published in transportation research journals such as IEEE annual conferences, Sustainability, Sensors, Transportation Research (Parts: B, C, and E), Journal of Advanced Transportation, Transportation Research Procedia, Transportation Science; European Journal of Operational Research, Cities, etc. Others include Expert Systems with Applications, Journal of Cleaner Production, Applied Energy, etc. Figure 3 refers to the network of journal bibliographic coupling.

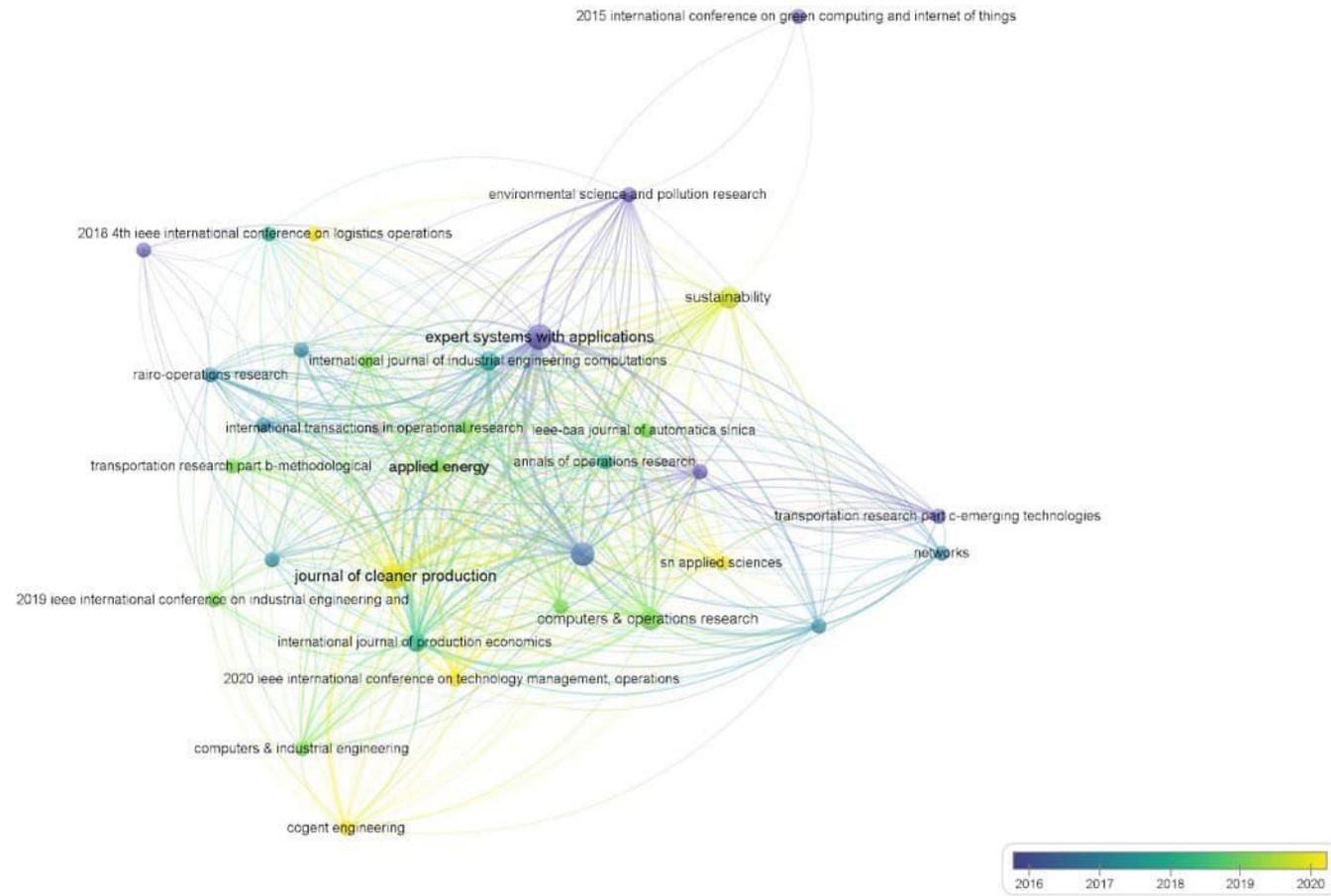

Figure 3: The network of journal bibliographic coupling

Contributing authors in GVRP publications were mainly from the countries shown in Figure 4. The U.S started investigating GVRP earlier than any other country as most early publications are from USA. China, as the first rank, has improved in GVRP research in recent years. Furthermore, five major European countries (Germany, England, Spain, Italy and France) have engaged considerably in the improvement of GVRP field. Iran, Canada and India have a high rate of publication, playing an influential role in developing GVRP solutions in the recent years. Probably, the high rate of publication in China, Europe and the U.S is associated with their share of the AFV market. As of December 2020, China had the largest stock of EVS, with 42% of the global plug-in passenger EV fleet in use. China also dominates the plug-in light commercial electric vehicle and electric bus deployment, with 65% of the global commercial EV fleet [3]. Europe had about 3 million plug-in passenger EVs by the end of 2020, accounting for 30% of the global stock [4]. It also has the second-largest electric light commercial vehicle fleet, with about 31% of the global stock in 2019 [3]. As of November 2021, the U.S has the third-largest share of the EV market, after China and Europe [5].

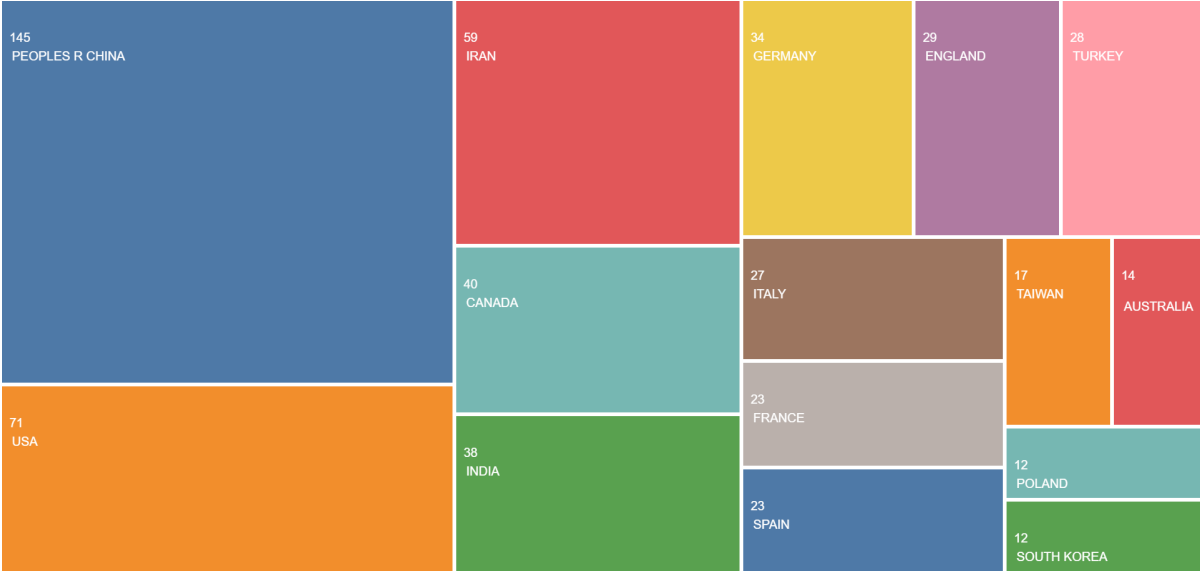

Figure 4: GVRP literature contribution by country

The terms used in the title, abstract and keyword lists of the green vehicle routing problem could be useful semantic lenses through the structure and composition of this literature and theme and nature of its studies. They could reflect the methods used in GVRP literature, the range of applications to which GVRP have been extended over various periods of time. Here, networks of keyword co-occurrence are analysed for GVRP publications as well as bursts of citations to keywords (or more specifically, to the articles where such keywords have been used). Two terms (which could be keywords or terms identified in the title or abstract) are considered co-occurred when they appear in the (keyword list or title/abstract of) the same publication. A map of term co-occurrence identifies those that have occurred most frequently and visualises those that have frequently co-occurred in closer spatial proximity [6]. Groups of terms that are highly related (i.e. those with strong co-occurrence relationship) form clusters of terms that could represent various sectors of the literature. A map of co-occurrence could be visualised in various formats including Figure 5, where the frequency of occurrence is presented with the years they have been appealing to scholars the most. This map indicates that the terms green logistics, vehicle routing problem and fuel consumption are amongst the recently occurred keywords in the publications. GA, ALNS, and TS were the most frequent algorithms in the GVRP literature.

Figure 5: Map of keyword co-occurrence for GVRP

# 3. GVRP Algorithms

In the literature, algorithms have been developed to resolve variants of Vehicle routing problems. The objective function of a VRP is defined based on the particular purpose of the study. Depending on applications, various types of GVRP are introduced, formulated and solved by exact or approximate methodologies. The main limitation of GVRP's proposed algorithms is generality. Specific heuristic and metaheuristic methods presented in one paper and one case study do not guarantee the effectiveness of solving other types of GVRP problems ; Therefore, in order to solve other variant problems more efficiently, more general methodologies are to be adopted.

**Exact methods** are currently applied in order to find an optimal solution for a few customers with limited capacity and fixed time windows [7, 8]. Popular exact procedures include direct tree search, dynamic programming, integer linear programming, etc. Although simple variants of real problems can be modelled with graphs naturally, even some simple variants of vehicle routing problems are considered to be NP-hard due to the size and frequency of real world VRPs in networks. So, the size of problems that can be solved, optimally, using mathematical programming or combinatorial optimization may be limited. Therefore, several recent studies have turned to approximate algorithms, evaluating their efficiencies in solving vehicle routing problems. Approximate algorithms can be heuristic, meta-heuristic, hybrid heuristic or machine learning methods.

**Heuristic methods** are used to resolve the vehicle routing problem according to the specific knowledge of the problem, which is in most cases suboptimal or close enough to a reliable solution [9]. In GVRP literature, heuristic methods can be split into constructive and improvement heuristics.

*Constructive heuristics* seek to propose an initial solution by providing either serial or parallel route construction [10]. Such solutions are constructed in a greedy way, which usually generates solutions slightly far from an optimal solution of the VRP. In this regard, the modified savings method is used to provide an initial solution of several types of GVRP, and especially Electric-VRP with the insertion of charging stations.

Traditional local search algorithms usually evaluate the whole neighborhood, but only perform one single move at each step. However, there are often many neighborhood moves in the current neighborhood that are independent of each other and can be simultaneously performed without interference [11]. The local search stops when no improvement in the solution can be noticed in the neighborhood of the incumbent solution, also named as local optima.

**Metaheuristic methods** can be defined as heuristics guiding other heuristics. These methods are either neighborhood-oriented (local) metaheuristics or population metaheuristics. *Neighborhood-oriented heuristics* keep exploring the neighborhood of the optimal solution. Simulated Annealing (SA), Tabu Search (TS), Variable Neighborhood Search (VNS), and Adaptive Large Neighborhood Search (ALNS) are among popular local metaheuristic methods [12]. *Population metaheuristic* methods are based on the natural selection procedure to evolve a population and let the fittest survive. Among them, Genetic Algorithms (GA), Ant Colony (AC), Bee Colony (BC), and Particle Swarm Optimization (PSO) are often used [13].

**Hybrid metaheuristic methods** take advantage of the meta-heuristic procedures to keep searching even after reaching the first local optima. In some cases, several metaheuristic and heuristic methods are combined and applied to a vehicle routing problem since using a specific approach leads to difficulties, such as a low-quality solution, trapping in local optima in search space, or high computation time. Therefore, several studies hybridize two or more algorithms to simultaneously employ strengths. Hybrid methods include exact-metaheuristic [14], metaheuristic-metaheuristic [15] and metaheuristic-heuristic [16] algorithms in order to obtain better results and add to the robustness of the solution. Moreover, several studies have attempted to solve GVRP with general exact solvers, such as CPLEX, Lingo and GAMS [17].

It should be noted that despite the development of exact methods, very few exact methods have been proposed for the EVRP and its extensions, which is a branch of GVRP. Exact methods are found to be inadequate to solve a large-scale optimization problem [18]. While researchers have used a population metaheuristic to solve the problem, only a few studies were able to generate high-quality solutions in a reasonable computational time. Most importantly, the use of emerging machine learning and data mining tools has been overlooked in the literature of GVRP algorithms.

Figure 6 refers to the overall percentage usage of each type of method in GVRP literature. It can be inferred that metaheuristic algorithms are increasingly recognized as a significant option, the most applicable methodology, and are currently receiving global attention from scholars and practitioners. Additionally, the small contribution of heuristic methods (about 7%) shows that most scholars tend to combine them into other methods due to local optima deficiency. This usually relates to these methodologies by concentrating on solving a specific problem.

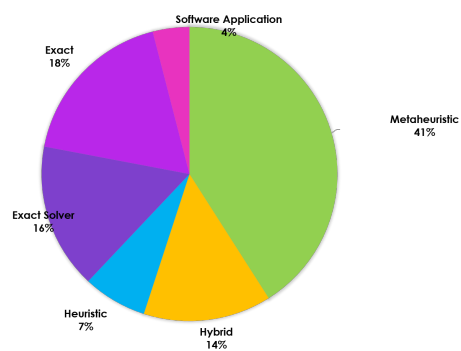

Figure 6: Usage of different solution methodologies in the existing literature

Figure 7 refers to the number of publications per year using each type of above-mentioned methodologies. This figure demonstrates the increasing trend in using metaheuristic algorithms and confirms the significance of this type of methodologies. Overall, the exact methods and metaheuristic algorithms are the methods most preferred by researchers. Furthermore, software applications and exact solvers have hit the bottom as the least popular methods due to their high computational time and complexity in the first years, while heuristic methods have become less popular, taking the place of exact solvers in the diagram.

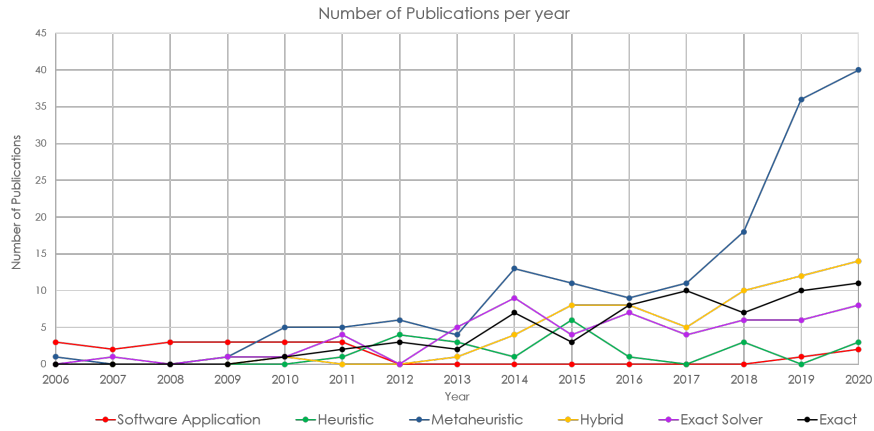

Figure 7: The number of publications per year based on methodology

## 4. GVRP Classification

We can divide GVRP literature into three main categories: (1) GVRP with conventional vehicles (2) GVRP with alternative fuel vehicles, and (3) GVRP with mixed fleet of vehicles. We also identified several variants as subcategories for the GVRP with Conventional Vehicles (CV) and Alternative Fuel Vehicles (AFV).

### 4.1. Conventional Fossil Fuel-Powered Vehicles Routing

The CO2 and NOx emissions problem has negative impact on the environment as well as on human health. In CVRP, any type of emission is considered in the objective function, with the main focus on minimizing the routing cost and polluting emissions. Several variants of GVRP with conventional vehicles have been explored by scholars in recent years. Figure 8 is associated with different GVRP with CV variants and their origination year in the literature. The figure reveals that while issues such as time dependency, refueling location and eco-driving have been investigated since the beginning of last decade, emerging issues like multi-objective, connectivity and especially automation has been overlooked and only in the last 2-3 years have the researchers found interest in examining GVRP considering connected and automated vehicles.

### 4.1.1. Multi-Objective Conventional GVRP

In this type of CVRP, more than 1 objective is taken into account. First introduced by Demir et al. [9], this variant aims to minimize both route cost and fuel consumption or vehicle kilometers travelled (VKT), or travel time or speed. In their study, they solved a bi-objective pollution routing problem (PRP), minimizing two conflicting factors: fuel consumption and driver time. Several scholars have recently developed Demir's methodology for solving PRP such as [19, 20]. Poonthalir and Nadarajan [19] resolved a multi-objective problem with varying speed constraints. Their model was able to minimize both route cost and fuel consumption, using particle swarm optimization with a new mutation operator called greedy mutation operator. Rauniyar et al. [20] developed a reliable solution methodology based on genetic algorithm to solve a bi-objective C-GVRP defined by Demir et al. Some studies have investigated multi-objectives like minimizing marginal cost and fuel consumption, VKT and travel time, combined with other variants. Alfaseeh et al. [21] and Djavadian et al. [22] incorporated connectivity and automation to their multi-objective Eco-routing study considering two different routing strategies as myopic and anticipatory routing. They showed that the anticipatory routing strategies outperform myopic ones due to the consideration of future traffic state in their routing calculations.

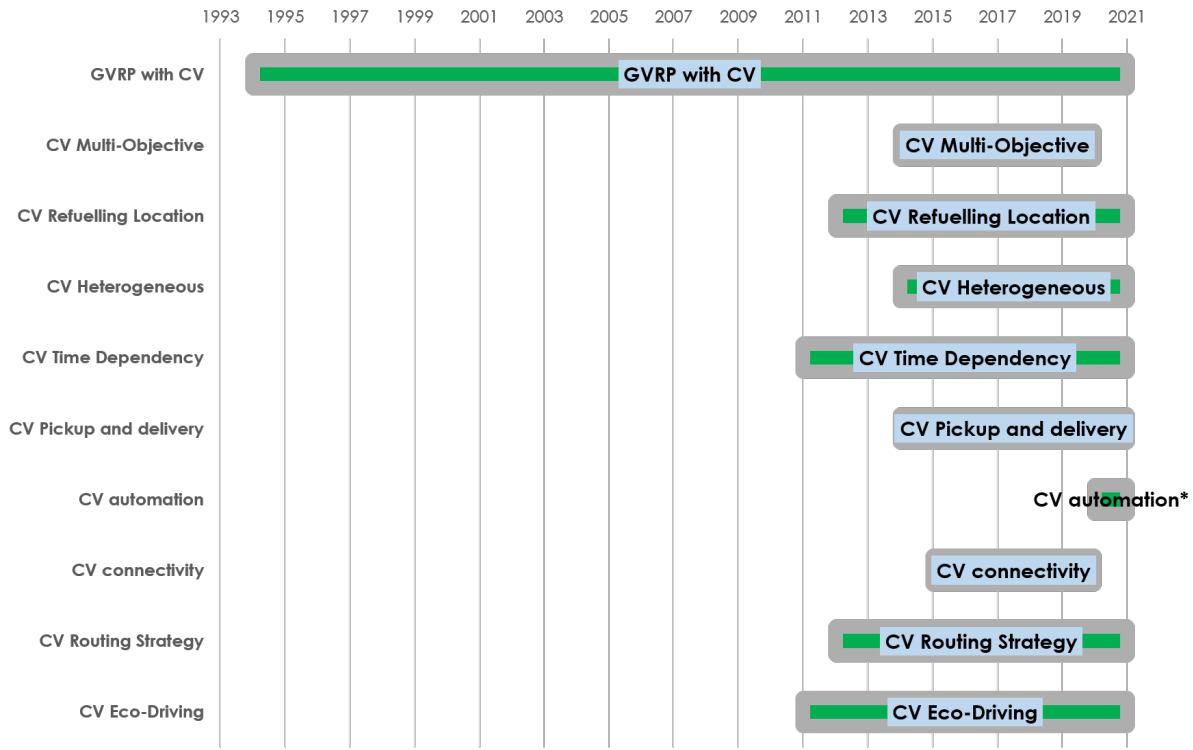

Figure 8: GVRP with CV variants origination

*4.1.2. Refueling Locations*

In several papers such as Yi et al. [23], constraints like road conditions, congestion, topography, vehicle load and their impact on route cost and fuel consumption have been included. Rezaei et al. [24] examined green vehicle routing problem with time windows constraint considering a heterogeneous fleet of vehicles and fuel stations. However, their methodology can be only applied to products with hard time windows such as like milk, meat, and newspaper.

*4.1.3. Heterogeneous Fleet*

This variant is associated with a fleet of vehicles with different capacities and costs, available for distribution activities. The problem is also known as the Mixed Fleet VRP or as the Heterogeneous Fleet VRP. Koc et al. [25] first modelled and solved a GVRP variant with a heterogeneous fleet. They applied a hybrid evolutionary metaheuristic method to solve the problem and concluded that in an urban area, the advantages of having a heterogeneous fleet outweigh those of a homogeneous one.

*4.1.4. Time Dependency*

Sometimes, GHG emission varies over time. Thus, emission can be a function of time because of travel speed variability creating time dependency in GVRP. This variant deals with the difficulties brought by the various factors of travel speed variability, such as the travel speed, congestion, and land use affecting $CO_2$ emissions [26]. Also, VRP with Time Windows (VRPTW) is to deliver the goods with time constraints and limited capacity of the vehicle fleet. This variant includes both soft/hard Time Windows and a combination of both. The complexity is that there are various uncertainties in this variant due to unexpected occurrences. Also, the algorithms used in the previous studies, were mostly applied to products that hard time windows are suitable for them such as meat, and newspaper. Franceschetti et al. [27] directly considered the effect of traffic congestion into conventional GVRP with hard time windows. Jabali et al. [28] investigated the two phases of free-flow traffic and congestion. They minimized the emissions per kilometer as a function of speed, developing a relationship between the reduction of emissions and marginal costs.

*4.1.5. Pickup and Delivery*

This is a more general issue in green logistics. The aim of this variant is to minimize operational delivery costs and the environmental impacts at the same time. In such cases, there are several issues to be considered. Some of them are: locating depots where vehicle with limited capacity is sent to deliver orders, minimizing emissions by scheduling customers. Tajik et al. [29] showed that considering fixed cost, the vehicle routes would better be composed of large vehicles, despite their higher emissions.

*4.1.6. Automation and Connectivity*

These variants are quite novel in the literature of GVRP. Issues regarding connectivity and Automation were only discussed by very few scholars such as [21, 22]. They aimed to minimize travel time, GHG and NOx emissions with different costing approaches and routing strategies. They asserted that as vehicular ad hoc network (VANET) impact on eco-routing is under the direct influence of ITS application, its sensitivity to the availability, robustness, accuracy, and temporal and spatial distribution of the network data should be further investigated. Furthermore, the authors believed that routing based on GHG as the objective offers a considerable reduction in average TT, average VKT, total GHG, and total NOx compared to the alternative where TT is the main objective.

*4.1.7. Routing Strategy*

Routing strategy can be an influential factor in optimization purposes. Inventory Routing allows suppliers to deliver their products to a given set of customers while optimizing inventory management, vehicle routing, and delivery schedule all at the same time. This is mostly applicable when delivering goods are subject to various constraints. It is assumed that emissions are associated with routing decisions, whereas waste is more linked to inventory decisions. As for future research direction, the multiproduct distribution, and non-deterministic consumption rates are among the features that have to be explored in order to minimize both the total cost of distribution and emission.

*4.1.8. Eco-Driving*

Eco-driving is defined as driving in such a way that minimizes emission by considering the dynamics of the traffic flow and safety measures. Zhou et al. [30] investigated fuel consumption models to evaluate eco-driving and eco-routing. They proposed that drivers often have more difficulty in applying eco-driving techniques on roads with high congestion. When traffic conditions are stable, eco-driving is more successful compared with facilities with lower speed limits and several roundabouts and ramps. Also, on higher speed limit roads, free-flow conditions can cause an increase in cruising speed, which lead to higher instant fuel consumption and consequently higher emission.

*4.2. Alternative Fuel Vehicle Routing*

The main goal of Alternative Fuel Vehicle Routing Problem (AFVRP) is to provide optimal routes with minimum energy consumption, time or cost for a fleet of alternative fuel vehicles, while considering their operation limitations such as limited driving distance and capacity. The number of publications on AFVRP is illustrated in Figure 9. There are more than 300 publications on AFVRP according to Scopus.

Alternative Fuel Vehicle Routing Problem (AFVRP) can be divided into six categories based on their fuel type, as illustrated by Figure 10. Figure 11 displays the frequency publications on AFVRP per year divided into 3 main categories of EVs, Hybrid EVs, and other alternative fuel types. Figure 12 illustrates the share of each type of fuel in AFVRP literature on the whole. Although there is a wide range of fuels associated with this type of GVRP, only a small fraction of the literature is focused on non-electric AFVRP. Also, from Figure 11, it can be inferred that 23 (18%), 90 (72%), and 12 (10%) studies belong to the AFVRP, EVRP, and Hybrid VRP (HVRP), respectively. Previous studies have mostly been focused on EVRP variants, and there exists a research gap on the other two variants of the problem, and specifically the HVRP variant.

Figure 13 shows the frequency of occurrence for various keywords, their co-occurrence, and the year they were most attractive to the researchers. This map indicates that the terms optimization, liquefied natural gas, refueling strategies, and time windows are amongst the most occurred keywords in the AFVRP publications. The term adaptive behavior is departed from the whole graph being used only once in the literature back in 2014.

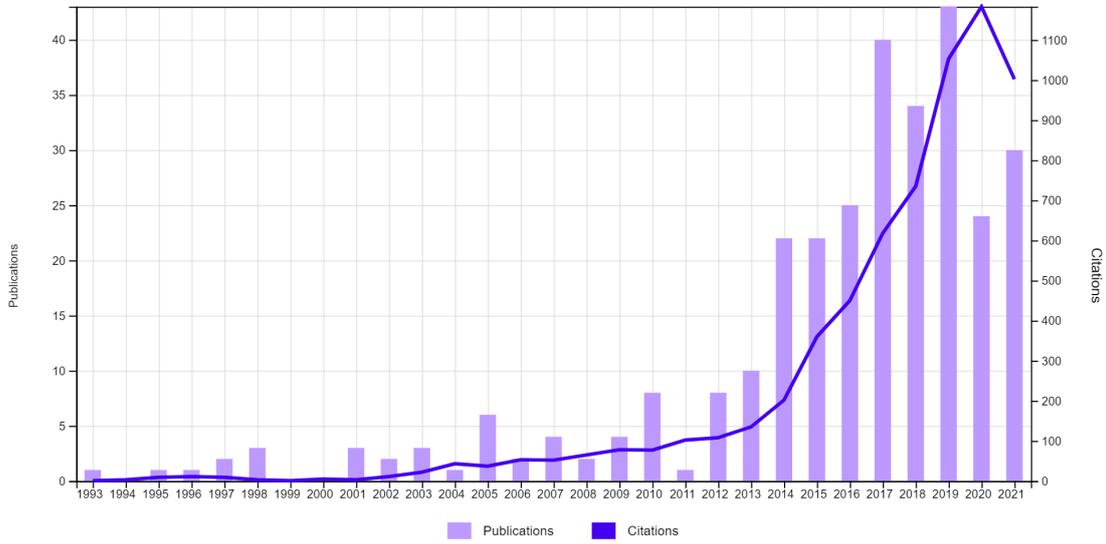

Figure 9: Alternative fuel vehicle routing publications per year

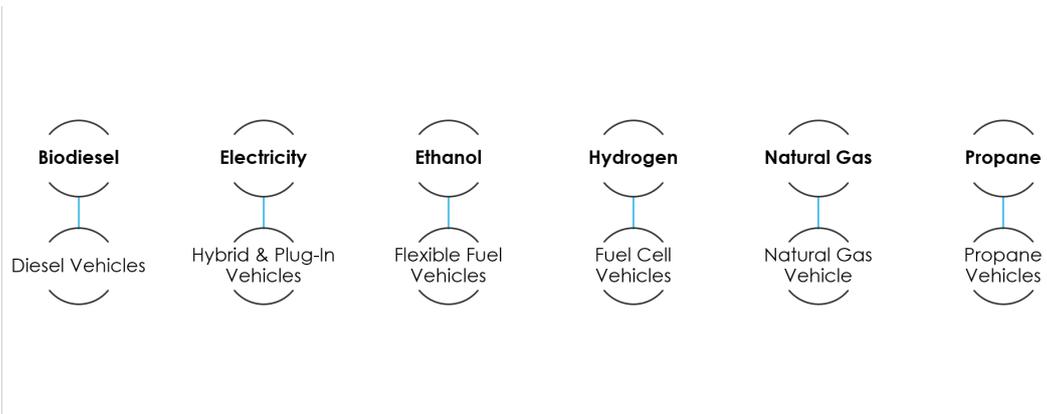

Figure 10: Alternative fuel vehicle classification

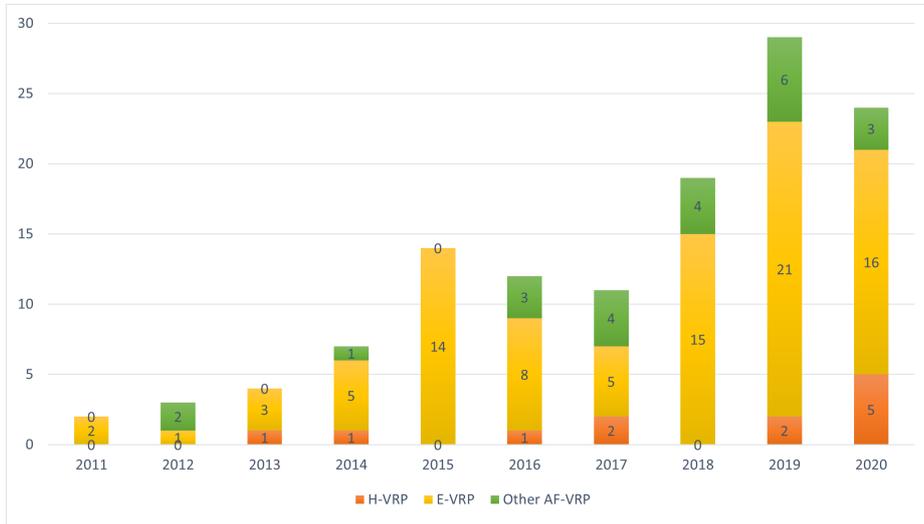

Figure 11: Number of EVs, Hybrid EVs, and other alternative fuel vehicles' publications per year

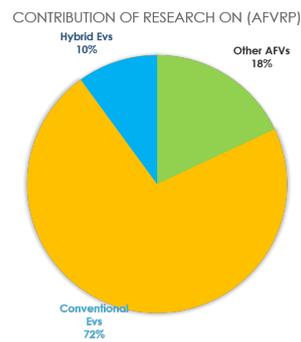

Figure 12: Percentage of AFVRP publications based on fuel type

Figure 13: Map of keyword co-occurrence for GVRP

Figure 14 is associated with different GVRP with AFV variants and their origination year in literature. The figure reveals that while issues related to refueling location have been investigated since the beginning of the last decade, emerging issues like partial recharging have only gained interest in recent years with the advancement of new charging technologies and battery swapping.

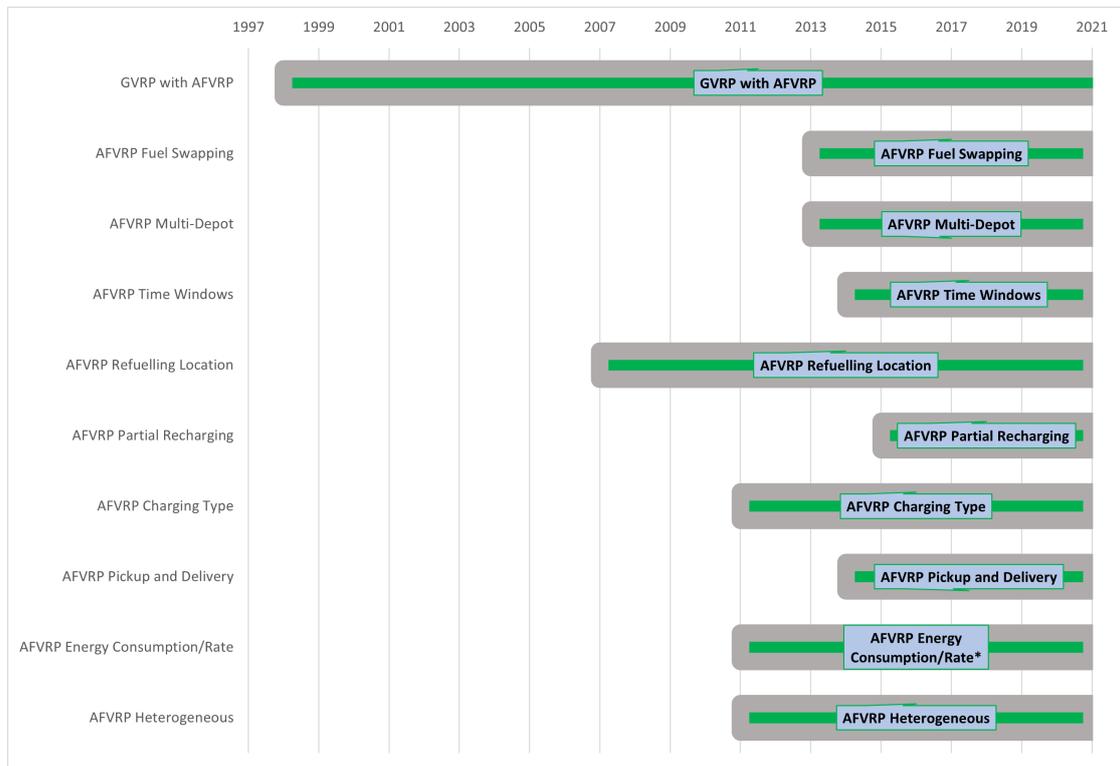

Figure 14: GVRP with AFVRP variants origination

*4.2.1. Time Windows*

This variant was considered earlier than other EVRPs, while it is more complex than the classical VRP with time windows. In this variant, a set of EVs deliver goods to customers, considering a given time window. Not visiting in pre-specified intervals decreases customers' satisfaction or may lead to infeasible solutions in practice. In one of the most recent papers associated with this variant, Keskin and Catay [31] considered stochastic waiting times at recharging stations with time windows. In another study, Schneider et al. [32] proposed a hybrid heuristic, which is a combination of a variable neighborhood search algorithm with a tabu search heuristic, considering limited vehicle freight capacities as well as customer time windows. Also, in a similar case study, Ding et al. [33] developed a heuristic method, according to variable neighborhood search and tabu search, by applying simple charging time adjustment processes to provide a more efficient solution. Moreover, Mao et al. [34] used an improved ant colony optimization (ACO) algorithm and hybridized it with insertion heuristic and enhanced local search to provide a solution for this problem, with respect to partial recharging and battery swapping. In addition, they proposed a new probabilistic selection model in ACO by considering the impact of both distances and time windows. In addition, Keskin and Catay [35] conducted more practical research by applying the full recharge restriction and allowing partial recharging (EVRPTW-PR). They defined this problem as a 0–1 mixed-integer linear program and proposed an Adaptive Large Neighborhood Search (ALNS) algorithm to solve it efficiently.

*4.2.2. Replenishment*

In this variant, it is usually assumed that a vehicle gets fully recharged every time it is at the refueling station–thus, continuing its service as long as its battery can support it. While EV battery packs are capable of supporting travel in the 100-mile range on a single charge, in order to replenish their batteries, they need to have access to charging stations. Meng and Ma [36] explored electric vehicle routing with soft time windows and offered two ways of power replenishment. The design of mobile charging stations have been investigated by Huang et al. [37]. In another study, Wen et al. [38] considered the problem of locating electronic replenishment stations for electric vehicles on a traffic network with flow-based demand to optimize the network performance. Baouche et al. [39] presented a new approach for the EV routing problem with recharging stage(s) along the way on the available charging stations to solve the autonomy limitation. Yet, issues like compatibility of battery EVs with chargers in their Charging Station(CS)s, and the impact of recharging at public or private CSs are still overlooked in the literature. Moreover, the time spent for recharging or battery swapping is a critical factor in this variant. A straightforward assumption usually used in vehicle routing studies is that the recharging time is constant across the whole network. Service time and charging time are also traditionally considered fixed.

*4.2.3. Charging Type*

In most of the existing EVRP papers, the battery charge level is a linear function of charging time, while in reality, this function is nonlinear. A practical linear charging estimation may lead to infeasible, unrealistic or expensive solutions. Therefore, Froger et al. [40] considered a realistic nonlinear relationship between the time spent on refueling and the amount of fuel consumed by the vehicle. Moreover, Xiao et al. [41] proposed a new model of the electric vehicle routing problem with respect to a general energy/electricity consumption function for EVs that factor in energy losses, nonlinear charging function with the piecewise linearization technique, efficient visits to charge stations, and continuous decision variables for speed, payload, travel time, recharging, etc. Karakatic [42] developed a Two-Layer Genetic Algorithm (TLGA) to overcome the capacitated Multi-Depot Vehicle Routing Problem by considering Time Windows (MDVRPTW) and Electric Vehicles (EV) with partial nonlinear recharging times. Zuo et al. [43] developed a new technical formulation for vehicle route selection and charging station visits, which reduces the formulation complexity. In addition, they proposed a new linearization method that applies a set of secant lines to surrogate the concave nonlinear charging function with linear constraints.

*4.2.4. Partial Recharging*

Conrad and Figliozzi [44] first considered the possibility of partial and overnight recharging for EVs at customer sites, which brought about one of the operational variants of the E-VRP. Subsequently, Felipe et al. [45] modelled and solved an EVRP with partial recharging for the first time. Since then, different aspects such as limited number of chargers in charging stations, various kinds of charging, fuzzy optimization models, etc. have been looked into. A method of charging is to establish some stations for swapping their batteries along the route either via stations or the road infrastructure. However, the impact of the different elements of this variant such as vehicle size, geographical configuration of site, recharge stations, autonomy, and recharging technologies have to be further checked. As one of the recent studies in this regard, Kancharla and Ramadurai [46] investigated capacitated and load-dependent discharging in a partial recharging EVRP. Also, several other combinations of this variant with time windows [47] and congestion tolls [48] have been studied to date. Hiermann et al. [49] proposed an electric vehicle routing problem combining conventional, plug-in hybrid, and electric vehicles.

*4.2.5. Location of Charging Station*

The driving range of AFVs is typically limited. Fuel cell electric vehicles are fueled with pure hydrogen gas stored in a tank on the vehicle and have a driving range of over 500 kilometres. Thus, the choice of routing with the focus on the driving range and long duration of the recharging process as well as the location of charging stations to utilize the necessary charging infrastructure is of significance. This variant accounts for the limited driving range of BEVs, which directly leads to the more frequent recharging needs at CS. However, only a few studies have considered the problem of CS capacity, due to the limited number of BEVs which can be charged at a CS, simultaneously. Only in recent years, have scholars started to combine a nonlinear charging process and CS location problem with distance constraint into the EVRP models [50]. They have proposed that waiting time at the charging stations can increase the total cost. Thus, considering waiting times, battery reservations, and adaptive routing with uncertainties in CS availability is an integral part of future EVRP research direction. Zhang et al. [51] also developed an Improved Whale Optimization Algorithm (IWOA), presented the locating problem of Electric Vehicle (EV) charging stations model with service risk constraints and applied IWOA to solve it. In another study, Bilal and Rizwan [52] investigated different approaches, objective functions, constraints and range of optimization techniques that were addressed by researchers for optimal placement of CS during recent years.

*4.2.6. Pickup and Delivery*

BEVs have a shorter driving range (160–240 km) compared to the driving range of conventional vehicles (480–650km). Therefore, they are more likely to be used on short distances or in urban areas where they are more effective than conventional vehicles due to their lower driving speed, lower noise production, and cost [53]. When the average route length is short, BEVs can be easily used as recharging can happen upon their return to the depot. BEVs are already being used by companies like DHL, and FedEx mostly for last-mile delivery of light goods as distances are short enough. Several studies have been conducted to investigate pickup and delivery effects in AFVs. For instance, Madankumar and Rajendran [54] presented two Mixed Integer Linear Programming (MILP) models for solving the Green Vehicle Routing Problems with Pickups and Deliveries in a Semiconductor Supply Chain. In another research, Ahmadi et al. [55] also investigated the importance of vehicle dynamics parameters in energy models for EV routing, especially in the Pickup-and-Delivery Problem (PDP). Grandinetti et al [56] developed a multi-objective mixed integer linear model for minimising the total travel distance, the total cost for the EVs used and the total penalty cost for the unsatisfied time windows.

*4.2.7. Energy Consumption Rate/Model*

Several energy consumption models were used by the literature, but only a few can predict realistic energy consumption at the road segment level, such as [57]. Additionally, Qi et al. [58] proposed a model to obtain an accurate link-level energy consumption estimation for EVs, considering the energy consumption under real-world traffic congestion on two proposed impact factors i.e. positive kinetic energy (PKE) and negative kinetic energy (NKE). Yi and Bauer [59] developed an adaptive multiresolution framework for electric vehicle (EV) energy consumption estimation with real-time capability. Accordingly, three key parameters, namely powertrain efficiency, wind speed, and rolling resistance, are adaptively estimated using a two-step nonlinear iterative algorithm. Basso et al. [60] introduced the Two-stage Electric Vehicle Routing Problem that incorporates improved energy consumption estimation with respect to detailed topography and speed profiles.

*4.2.8. Multi-Depot*

This variant arises as a practical and functional issue in VRP, where vehicles are dispatched from and returned to one of the multiple depot locations. Therefore, apart from the routing choice, it is necessary to decide from which depot the goods are going to be delivered. Only 4% of the publications have incorporated this variant, including [61]. In this regard, Zhang et al. [62] used an Ant Colony System-based metaheuristic to find a solution to A Multi-depot Green Vehicle Routing Problem (MDGVRP). Schneider et al. [63] considered the Multi-Depot Electric Vehicle Location Routing Problem with Time Windows (MDVLRP).

*4.2.9. Fuel Swapping*

This variant is defined to establish some stations for swapping batteries and other fuel types which correspond to increased recharging speed and reduced time loss. A battery swapping operation is faster than recharging operation, by taking only 10 minutes. Also, the used-up batteries can be recharged at night when electricity is charged at a discount. This variant accounts for over 15% of the publications on AFVRP including [64]. Some scholars have proposed (recharging) or battery swapping services can be available at all or some of the customers' sites as well [65]. In other studies, such as [66], the authors developed robust optimization models that aid the planning process for deploying battery-swapping infrastructure. Also, Wu [67] reviewed the state-of-the-art battery swapping station (BSS) literature and business models, where the BSS offers a recharged battery to an incoming EV with a low state-of-charge. First, four operation modes are presented i.e., a single BSS, multiple BSSs, an integrated BSS and battery charging station (BCS), and multiple BSSs and BCSs. BSS problems in routing are surveyed in different operational areas including charging schedule, construction and planning, dispatching and routing optimization, and power management.

*4.2.10. Upstream Effect*

Both the EU and U.S.A boasted about electric vehicles producing zero emission. However, unlike conventional vehicles, a significant proportion of the emissions produced by electric vehicles occurs 'upstream', i.e., when the electricity is produced at the source. Thus, current regulations, which only account for exhaust emissions, do not fully capture the GHG emissions from EVs. If retaining the 0g/km rate is an actual goal, it is required that manufacturers buy carbon credits to compensate for the 0g/km rate or switch to a full life cycle analysis. Lutsey and Sperling [68] found that if upstream emission's effect are considered, an EV powered from the American electricity grid produced an average of 56% less CO2 emission than a similar brand new petrol car (62g/km compared with 142g/km). Sen et al. [69] uses a hybrid life-cycle assessment method to analyze and compare alternative fuel-powered heavy-duty trucks (HDTs) applying a Monte Carlo simulation to account for the uncertainty in the data.

## 5. New Insights and Perspectives into GVRP

Even considering a pervasive application of alternative fuel vehicles and their role in emission reduction, a large number of vehicles remain on the roads, producing congestion and polluting emissions. Hence, to ultimately reduce the emission, we must consider emerging solutions that can also take into account the increasing customer demand for vehicles, while considering environmental aspects.

*5.1. Technological Consideration*

Several emerging technologies are opening up new research directions in the context of GVRP. Here we discuss the most promising ones.

*5.1.1. Drones*

Unmanned Aerial/Ground Vehicles (UAV/UGV) or drones are an emerging technology solution to the last-mile delivery problem. Drones are either controlled by a remote controller or an on-board intelligence. They have the potential to cut the pollution caused by trucks on a congested road network by utilizing the unused airspace. A mixed fleet of drones and trucks can also be used for delivery as proposed by Wang el al. [70]. Drone technology is considered more reliable and faster since neither are drones affected by road congestion nor traffic accidents on the networks. Several research directions can be considered. First, the environment that the UAV/UGV is operating in should be specified. That is, within which obstacles and urban air mobility limitations the vehicle is making its way. To do so, vehicle flight planning and optimization are required and considered as a path planning problem. As a result, yet another version of GVRP can be formulated where the drone must optimize the path, while conserving the fuel and minimizing emissions. Sarath et al. [71] discussed several techniques to achieve UAV path detection, planning and obstacle avoidance for real-time communicative environments. Alongside obstacle avoidance, there are some aerial restrictions in the path planning of drones, such as no-fly zone areas. Feng [72] proposed an improved method to achieve path planning for UAVs in complex surroundings. However, the fuel consumption and emissions dimension is rarely considered in the path planning problem. Secondly, it is noted that the wind and weather play a critical role in the flight planning of the drones, especially if they are small and light. Cheng et al. [73] introduced a distributionally robust optimization model to solve a two-period drone scheduling problem with uncertain flight times, which can be implemented in a data-driven framework using historical weather information. The operation of drones for the last-mile delivery, where the mobile base-station is a truck, creates another emerging GVRP, where the integrated truck and drone operations need to be optimized in a dynamic environment. As the drone traffic is expected to increase in the near future, we expect that cities will regulate the airspace more and may consider extending the existing 2D on-land road network to 3D land and air road network. This will result in further potential research directions in GVRP.

*5.1.2. Distributed Systems*

In recent years, distributed ledger technologies, for instance, blockchain have been used in transportation and logistics to manage information. Lopez and Farooq [74] proposed a blockchain based smart mobility data market (BSMD) that provides the underlying framework for the use of distributed ledgers in transportation applications. Eckert et al. [75] developed a carbon credit (C2) market for multimodal passenger mobility using BSMD, where individuals could track and trade C2s based on their mode, trips, and availability of credits. Such a market has the potential to be used in GVRP where the emissions are minimized not only based on the cost objective, but also with cooperation/competition among the individuals in the market by exchanging carbon credits with dynamic pricing. Due to the distributed and dynamic characteristics of the Internet of Vehicles (IoV), content-centric decentralized vehicular named data networking (VNDN) has become more suitable for content-oriented applications in IoV [76]. The existing centralized architecture is prone to the failure of single points, which results in trust problems in key verification between cross-domain nodes, consuming more power and reducing the lifetime. Focusing on secure key management and power-efficient routing, [76] proposed a blockchain-based key management and green routing scheme for VNDN.

Due to advancements in information and communication technologies, there is a strong focus on developing highly intelligent intersections in urban areas that can control and route traffic. Farooq and Djavadian [77] proposed a distributed traffic management system where the intersections actively cooperate and exchange information among each other via Infrastructure to Infrastructure (I2I) communication technology. These intelligent intersections (I2s) use the information to predict traffic conditions and proactively route vehicles from origin to destination via vehicle to infrastructure (V2I) communication technology. We are of the view that such distributed and intelligent systems can be used to develop new solutions to GVRP.

*5.2. Fuel Technologies*

Advanced research on lithium-ion batteries done over several decades has resulted in high energy density, high cycle life, and high-efficiency batteries. However, the research is still ongoing on new electrode materials to enhance the performance of energy density, power density, cycle life, safety, and cost. The current generation of anode and cathode materials are suffering from several issues, including, slow Li-ion transport, high volume expansion, limited electrical conductivity, low thermal stability, dissolution or other unfavorable interactions with electrolyte, and mechanical brittleness [78]. Several approaches have been developed to solve these problems. A variety of intercalation cathodes have been available on the market, and conversion material technology is going to become more common. In terms of GVRP, lithium-ion battery electrode materials' technological advances would result in solutions with less stringent constraints in terms of trip lengths and better utilization of the capacity of the vehicle, especially when the demand is highly stochastic.

In the view of the Fraunhofer Institute, synthetic fuels and drive technologies such as hydrogen combined with the fuel cell has the potential to play a crucial role in the future of transportation. It is expected that such a role might be negligible in private vehicles, but significant in long-distance and heavy-duty vehicles used for goods movement. However, the drawbacks of hydrogen-based fuel cell technology should not be overlooked as it is very costly in terms of efficiency and operating costs. Horváth et al. [79] study, compares two types of EVs from the customer's point of view. In their study they had a detailed investigation carried out into whether battery- or hydrogen-powered electric vehicles will become ubiquitous in the future. The question of which energy storage system has the best efficiency and is the most cost-effective one for powering electric vehicles is still unanswered. With BEVs, only eight percent of the energy is lost upstream, and another 18 percent is wasted to convert the electrical energy to drive power. Depending on the model, BEVs' efficiency is about 70 to 80 percent. In the case of the hydrogen-powered EVs, 45 percent of the energy is already lost during the production of hydrogen through electrolysis. Of this remaining original energy, another 55 percent is lost to convert hydrogen to electricity. This means that the hydrogen-powered EVs' efficiency is about 25 to 35 percent, depending on the model. The efficiency is even worse with alternative fuels. The efficiency in this case is only 10 to 20 percent, which can convey the meaning that the use of Hydrogen would therefore be a mistake for passenger cars [79]. Therefore, Horváth et al. implied that investments should rather focus on long-distance and heavy duty transportation where ecological and economic constraints play an important role.

Although alternative fuels may become prevalent in the long term, yet it is unlikely that they would completely substitute the fossil fuels in the near future. Given the scarcity of such resources, alternative fuels should be prioritized for different transportation sectors to which it is cost-effective. This is not only because alternative fuels are competitive in those sectors but also because it is hard to decarbonize them. On the other hand, handling hydrogen from storage to transportation is difficult as it requires additional infrastructure such as hydrogen grid, and additional transformation on the demand side like fuel cells for heavy-duty transportation. All in all, considering the complexity of the adoption of these fuel technologies and the scarcity of vehicles that are powered from such technologies, green vehicle routing problem is less likely to be affected by studies in AFV areas rather than EVs and HEVs in the near future. Given that to inspect different aspects of GVRP with fuel cell vehicles, a sufficient rate of them is required on transportation network. To the best of our knowledge, very few studies have looked into this area of AFVRP. In addition, since the use of alternative fuel in heavy duty transportation is preferred to passenger vehicles, it is observed that scholars tend to investigate AFVRP in urban transit, logistics and air and rail transportation. Therefore, AFVRP studies will grow vastly based on the demand and penetration rates of alternative fuel vehicles.

*5.3. Methodological Considerations*

The application of machine learning methods in green vehicle routing and optimization has attracted scholars' attention in the last few years. While demonstrations using AI in GVRP are rare, in a recent publication, Guiladi and Eriksson [80] referred to the dynamic routing of electric commercial vehicles with a large amount of data considering random customer requests when predicting the optimal route. This work introduced artificial intelligence applied to routing and energy prediction of electric vehicles with a Deep Q-Learning method proposed to solve the problem. In a similar study, Chen et al. [81] proposed a Deep Q-Learning method to assign customers to vehicles and drones for same-day delivery. Although reinforcement learning in GVRP is not commonly used, it is considered a powerful tool for considering generalized GVRP studies considering different variants and applications of it where specific heuristic and metaheuristic methods can't be generalized to another problem. Also, it is a powerful method of dealing with uncertainties in the real world GVRP problems. In a study by Basso [82], it has been shown that the reinforcement learning (RL) method could save on average 4.8% (up to 12%) energy by planning the route and charging anticipatively, rather than the deterministic online reoptimization method. The research addresses the dynamic Stochastic Electric Vehicle Routing Problem (DS-EVRP) with a safe reinforcement learning method to solve the problem. In similar studies, a deep reinforcement learning method is proposed to minimize the total cost of a fleet of electric-autonomous vehicles for ride hailing services [83], where the complete system state is approximated using neural networks. Recently, approaches like chaos theory, quantum computation, and fuzzy logic have been introduced in GVRP, which can help deal with uncertainties. From the fuzzy theory perspective, most studies in this area consider the fuzzy chance-constrained mixed integer non-linear programming model ([84], [85]). In terms of stochastic optimizations, however, only a few studies have shown an application of this approach in the green freight transportation context [86].

Diversification is another crucial factor to the performance of the population-based algorithm, but the initial population in the many methods such as bee colony is generated using a greedy heuristic, which has insufficient diversification. Therefore the ways in which the sequential optimization for the initial population drives the population toward improved solutions are examined [87].

*5.4. Emerging Services and their Requirements*

Cars travelling on the roads are not completely exploited and can be used to deliver goods as well as passengers. Crowd shipping is replacing conventional delivery companies with occasional drivers using their personal vehicles to deliver goods. Archetti et al. [88] put forward the concept of VRP with occasional drivers (VRPOD). Macrina et al. [89] added time windows, multiple deliveries for origin-destination pairs, and split and delivery policy in several publications [90, 91]. It also introduced a VRP with a mixed fleet of CVs and EVs [89]. This strategy has proven to positively impact emission and routing costs reduction as well as offering a higher reliability and customer satisfaction level. Although crowd shipping companies have expanded their business in the last several years. However, those companies mainly have businesses within some metropolitan areas. A system may work for the last-mile delivery, but how it will perform for inter-city delivery is still to be investigated [92]. Another concern is that a business model may have differential performance in different contexts, possibly due to cultural differences and infrastructure networks that support CS markets [93]. In that a case, the promising application areas challenge stakeholders on both supply and demand sides (e.g., market segments, network issues), operations and management (e.g., reverse logistics) to implement CS systems that function collaboratively, dynamically, and sustainably.

## 6. Conclusion

To the best of our knowledge, no existing study has comprehensively considered all GVRP variants and its future direction. Review studies are crucial for understanding the existing body of literature and, in this regard, a clear gap exists. Other operational constraints driven from GVRP variants such as site-dependent GVRP, and periodic GVRP are to be investigated in the future. Also, if retaining the zero emission rate is a goal in reality, it is required that manufacturers buy carbon credits to make up for their emission, or switch to a full life-cycle analysis. However, this requires more data collection on GHG emission from the electricity grid, which may lessen the popularity of electric vehicles among manufacturers, leading to a sharp fall in their investment in EV market. Also, it can be concluded that the main focus of previous literature is on the operational level routing decision, overlooking other aspects associated with the supply chain management such as network design, road tolls, reliability index, etc. While there is a large body of literature developed on EVRP in recent years, AFVRP is yet to be further explored – especially in terms of mixed fleet, connectivity and autonomous vehicles. The other limitation of AFVRP is that some alternative fuels are not cost-effective for private vehicles. Thus, the demand side of AFVRP is yet to be developed based on fuel technologies.

The forecasting methods such as data-driven machine learning and reinforcement learning methods were rarely investigated in GVRP. However, metaheuristic methods are dominant for all variants of GVRPs due to their less time consumption and reasonable precision, indicating the significance of this approach compared with others. Nevertheless, there is no consensus among scholars about the reliability of such methodology. Also, exact methods are rarely applied in time dependent VRP. Opportunities in solution methodologies such as the application of Deep Q-Learning method, quantum computing, chaos theory, reinforcement learning, etc. are still not popular. Further GVRP research using such methods is highly recommended due to their strength and robustness in forecasting different parameters. In addition, the advancement in information and communication technologies has given rise to distributed systems related issues which still need further development in terms of Internet of Vehicles (IoV) and transportation infrastructures.

Opportunities for considering uncertainties are provided as previous studies considered time-dependent concepts to deal with uncertainties, while demand, speed, and travel time are still the most important parameters remained as uncertainty in this field. In this regard it can be noticed that other non-deterministic parameters such as uncertain travel time have been neglected. Sustainability related indices such as social concerns, customer willingness, driver pattern, and operational risks are among the uncertainties of the supply chain network, which are totally overlooked by the former research. In addition, to develop a more promising routing system, especially on a macroscopic scale with plenty of complex and hazardous areas, a road network should be designed to enhance network reliability and overcome vehicle routing issues. There are various issues, such as weather conditions, restricted zones areas, and several moving obstacles, that must be taken into account. In other words, the more real-world challenges to be considered, the more operational routing network will be achieved. Opportunities for multiple-objective approaches still exist as in reality, there is rarely a single-objective problem, while only 20% of the publications consider multi-objective problems, which confirms the need for multi-objective optimization problems in the green transportation context. Last but not least, new opportunities in urban logistics have paved the way for UAV/UGVs to influence GVRP in the last mile delivery. Crowd shipping in all the transportation modes has also been another strategic element of green logistics in recent years, saving the delivery time and cost, while producing less GHG emission. It is expected that in the near future, GVRP will not be limited to the land roads networks, but will also involve virtual air networks.